\documentclass[11pt]{amsart}

\usepackage{amscd,verbatim}
\usepackage{amsmath, amssymb,url}
\usepackage{amsfonts}
\newcommand{\de}{\partial}

\newcommand{\ddbar}{i \partial \overline{\partial}}
\newcommand{\Ric}{\mathrm{Ric}}
\newcommand{\ov}[1]{\overline{#1}}

\newcommand{\ti}[1]{\tilde{#1}}
\newcommand{\vp}{\varphi}
\newcommand{\vol}{\mathrm{Vol}}

\newcommand{\ve}{\varepsilon}

\renewcommand{\leq}{\leqslant}
\renewcommand{\geq}{\geqslant}

\newcommand{\be}{\begin{equation}}
\newcommand{\ee}{\end{equation}}

\begin{document}
\newcounter{remark}
\newcounter{theor}
\setcounter{remark}{0}
\setcounter{theor}{1}
\newtheorem{claim}{Claim}
\newtheorem{theorem}{Theorem}[section]
\newtheorem{lemma}[theorem]{Lemma}
\newtheorem{corollary}[theorem]{Corollary}
\newtheorem{conjecture}[theorem]{Conjecture}
\newtheorem{proposition}[theorem]{Proposition}
\newtheorem{question}[theorem]{Question}
\newtheorem{defn}{Definition}[theor]
\theoremstyle{definition}
\newtheorem{rmk}{Remark}[section]

\newenvironment{example}[1][Example]{\addtocounter{remark}{1} \begin{trivlist}
\item[\hskip
\labelsep {\bfseries #1  \thesection.\theremark}]}{\end{trivlist}}

\title{Immortal solutions of the K\"ahler-Ricci flow}
\author{Valentino Tosatti}
\address{Courant Institute of Mathematical Sciences, New York University, 251 Mercer St, New York, NY 10012}
\email{tosatti@cims.nyu.edu}
\dedicatory{In memory of Steve Zelditch}
\begin{abstract}
We survey some recent developments on solutions of the K\"ahler-Ricci flow on compact K\"ahler manifolds which exist for all positive times.
\end{abstract}

\maketitle

\section{Introduction}
Let $(X^n,\omega_0)$ be a compact K\"ahler manifold, and consider the K\"ahler-Ricci flow
\begin{equation}\label{krf}
\frac{\de}{\de t}\omega(t)=-\Ric(\omega(t)),\quad\omega(0)=\omega_0.
\end{equation}
Up to rescaling the time parameter by a factor of $2$, this flow is precisely Hamilton's Ricci flow \cite{Ha} starting at the K\"ahler metric $\omega_0$. As noticed already by Hamilton \cite[p.257]{Ha}, the evolved metrics remain K\"ahler (with respect to the fixed complex structure on $X$), and so we can identify them with their K\"ahler forms $\omega(t)$.

A key observation of Cao \cite{Ca} is that K\"ahler-Ricci flow, which is an evolution equation for a metric tensor, is in fact equivalent to a parabolic PDE for a scalar function. Indeed, if we let
\begin{equation}
\alpha(t):=\omega_0-t\Ric(\omega_0),
\end{equation}
which is a family of closed real $(1,1)$-forms (positive definite for $t>0$ sufficiently small, but not necessarily for $t$ large), then it is easy to see (see e.g. \cite[\S 3.2]{To}) that a family $\omega(t)$ of K\"ahler metrics on $X$ (with $t\in [0,T)$ for some $T>0$) solves \eqref{krf} if and only if it is of the form
\begin{equation}
\omega(t)=\alpha(t)+\ddbar\vp(t),
\end{equation}
where the smooth functions $\vp(t)$ (which vary smoothly in $t\in [0,T)$) satisfy
\begin{equation}\label{ma1}
\left\{
                \begin{aligned}
                  &\frac{\de}{\de t}\vp(t)=\log\frac{(\alpha(t)+\ddbar\vp(t))^n}{\omega_0^n}\\
                  &\vp(0)=0\\
                  &\alpha(t)+\ddbar\vp(t)>0.
                \end{aligned}
              \right.
\end{equation}
This is a (strictly) parabolic complex Monge-Amp\`ere equation, which is a fully nonlinear PDE (as long as $n\geq 2$). Thus, short time existence of a unique smooth solution $\vp(t)$ is guaranteed by standard PDE theory, with $t\in [0,T_{\rm max})$ where $0<T_{\rm max}\leq+\infty$ is the maximal time of (forward) existence (see also Hamilton's paper \cite{Ha} for the original approach, which works for general Ricci flows on closed manifolds). As an aside, observe that one can also replace $\alpha(t)$ by any other smoothly varying cohomologous closed real $(1,1)$-form, and obtain another equivalent form for the complex Monge-Amp\`ere equation \eqref{ma1}. Sometimes a more judicious choice of $\alpha(t)$ is preferable in order to obtain estimates for $\omega(t)$, see e.g the author's lecture notes \cite{To}.

A solution $\omega(t)$ of the K\"ahler-Ricci flow \eqref{krf} is called {\em immortal} if $T_{\rm max}=+\infty$. The goal of this note is to give an overview of what is known about such immortal solutions, and formulate some well-known open problems.

\subsection*{Acknowledgments} We are grateful to Man-Chun Lee and to the referee for comments on an earlier draft. Part of this work was carried out during a visit to the Center of Mathematical Sciences and Applications at Harvard University, which we would like to thank for the hospitality. The author was partially supported by NSF grants DMS-2231783 and DMS-2404599. It is with immense sadness that we mourn our dear friend Steve Zelditch, whose boundless enthusiasm for mathematics and life has deeply touched all those who knew him.

\section{When do immortal solutions exist?}
Suppose $(X^n,\omega_0)$ is a compact K\"ahler manifold and $\omega(t),t\in [0,T_{\rm max}),$ is a solution of the K\"ahler-Ricci flow \eqref{krf}, with maximal existence time $0<T_{\rm max}\leq+\infty$. The Ricci form $\Ric(\omega)$ of any K\"ahler metric $\omega$ on $X$ is a closed real $(1,1)$-form whose cohomology class
\begin{equation}
  [\Ric(\omega)]\in H^{1,1}(X,\mathbb{R}),
\end{equation}
is independent of $\omega$, and equals $2\pi c_1(X)$, where $c_1(X)$ denotes the first Chern class of the manifold $X$. Thus, if in the flow equation \eqref{krf} we pass to cohomology classes, we get the ODE
\begin{equation}\label{krf2}
\frac{\de}{\de t}[\omega(t)]=-2\pi c_1(X),\quad[\omega(0)]=[\omega_0],
\end{equation}
whose solution is given by
\begin{equation}
[\omega(t)]=[\omega_0]-2\pi tc_1(X).
\end{equation}
In particular, for $t\in [0,T_{\rm max})$ the cohomology class $[\omega_0]-2\pi tc_1(X)$ contains a K\"ahler metric, and thus lies in the K\"ahler cone
\[\mathcal{C}_X\subset H^{1,1}(X,\mathbb{R}),\]
of cohomology classes which contain a K\"ahler metric. The following theorem of Tian-Zhang \cite{TiZ}, which improved earlier results of Cao \cite{Ca} and Tsuji \cite{Ts,Ts2}, gives a nice cohomological characterization of $T_{\rm max}$ (see also \cite[Theorem 3.3.1]{SWL}, \cite[Theorem 3.1]{To} and \cite[Theorem 3.1]{We} for detailed expositions):

\begin{theorem}\label{tz}
Given $(X^n,\omega_0)$ a compact K\"ahler manifold, the maximal existence time $T_{\rm max}$ of the K\"ahler-Ricci flow \eqref{krf} is given by
\begin{equation}
T_{\rm max}=\sup\{t>0\ |\ [\omega_0]-2\pi tc_1(X)\in\mathcal{C}_X\}.
\end{equation}
It follows that $T_{\rm max}=+\infty$ is equivalent to $-c_1(X)\in\ov{\mathcal{C}_X}$.
\end{theorem}

This shows that if $X$ admits a K\"ahler metric $\omega_0$ for which the K\"ahler-Ricci flow starting at $\omega_0$ is immortal, the same will hold for {\em any} other K\"ahler metric on $X$.
Now, cohomology classes in the closure $\ov{\mathcal{C}_X}$ are traditionally called {\em nef}, and $-2\pi c_1(X)=c_1(K_X)$ where $K_X=\Lambda^n T^{1,0}X^*$ denotes the canonical bundle of $X$. Thus, the condition $-c_1(X)\in\ov{\mathcal{C}_X}$ is often stated by saying that the canonical bundle $K_X$ is nef. When $X$ is projective algebraic, such manifolds are also known as smooth minimal models in birational geometry.

In summary, compact K\"ahler manifolds which support immortal solutions of the K\"ahler-Ricci flow are exactly those for which the canonical bundle is nef, and on such manifolds all solutions of the K\"ahler-Ricci flow are immortal. Examples of such manifolds include compact Riemann surfaces of genus $g\geq 1$, smooth hypersurfaces in $\mathbb{P}^{n+1}$ of degree $\geq n+2$, and products of these.

Using the Calabi-Yau Theorem \cite{Ya} we can give a more geometric reformulation of the condition that $K_X$ be nef. Indeed, after fixing an arbitrary K\"ahler metric $\omega$ on $X$, it is easy to see \cite[Lemma 2.2]{To} that $-2\pi c_1(X)$ is nef if and only if for every $\ve>0$ there is a smooth closed real $(1,1)$-form $\eta_\ve$ cohomologous to $-2\pi c_1(X)$ such that $\eta_\ve\geq -\ve\omega$ holds on $X$. By the Calabi-Yau Theorem \cite{Ya} we can find a K\"ahler metric $\omega_\ve$ such that $-\Ric(\omega_\ve)=\eta_\ve\geq -\ve\omega,$ or in other words
\begin{equation}\label{eps}
\Ric(\omega_\ve)\leq \ve\omega.
\end{equation}
In other words, $K_X$ is nef if and only if for every $\ve>0$ there is a K\"ahler metric $\omega_\ve$ on $X$ that satisfies \eqref{eps} (i.e. $\omega_\ve$ has ``almost nonpositive Ricci curvature'' in a certain sense). The choice of $\omega$ is of course irrelevant, as any two K\"ahler metrics on $X$ are uniformly equivalent.

Combining Theorem \ref{tz} with Demailly-P\u{a}un's numerical characterization of $\mathcal{C}_X$ in \cite{DP}, we obtain the following geometric characterization of $T_{\rm max}$, which was conjectured by Feldman-Ilmanen-Knopf in \cite[p.204]{FIK}:

\begin{theorem}\label{tz2}
Given $(X^n,\omega_0)$ a compact K\"ahler manifold, the maximal existence time $T_{\rm max}$ of the K\"ahler-Ricci flow \eqref{krf} is given by
\begin{equation}
T_{\rm max}=\inf\left\{t>0\ \bigg|\ \text{there exists }V\subset X \text{s.t.}\int_V([\omega_0]-2\pi tc_1(X))^{\dim V}=0\right\},
\end{equation}
where here $V$ is a closed irreducible positive-dimensional analytic subvariety of $X$.
\end{theorem}
In other words, if the flow ceases to exist at a finite time $T_{\rm max}$, there has to be such a subvariety $V$ whose volume in the induced evolving metrics
$$\mathrm{Vol}(V,\omega(t))=\frac{1}{(\dim V)!}\int_V\omega(t)^{\dim V}=\frac{1}{(\dim V)!}\int_V([\omega_0]-2\pi tc_1(X))^{\dim V}$$
is converging to zero as $t\uparrow T_{\rm max}$. The union of all such subvarieties is itself a closed analytic subvariety, and it is in fact equal to the singularity formation set of the flow (where the curvature of $\omega(t)$ is blowing up), as proved by Collins and the author \cite{CT} solving another conjecture in \cite{FIK}.

\section{Volume growth of  immortal solutions}
Given $(X^n,\omega_0)$ a compact K\"ahler manifold with $K_X$ nef, we know from Theorem \ref{tz} that the solution $\omega(t)$ of the K\"ahler-Ricci flow \eqref{krf} starting at $\omega_0$ is immortal. The main (vague) question is then the following:

\begin{question}\label{mq}
What is the behavior of $(X,\omega(t))$ as $t\to+\infty$?
\end{question}

To start getting a feeling for this question, let us examine the behavior of the total volume of $X$ with respect to the evolving metrics. Using Stokes' Theorem we have
\begin{equation}\label{v1}
\vol(X,\omega(t))=\frac{1}{n!}\int_X\omega(t)^n=\frac{1}{n!}\int_X([\omega_0]-2\pi tc_1(X))^n,
\end{equation}
and this is a polynomial in $t$ of degree equal to the {\em numerical dimension} of $K_X$
\begin{equation}
m:=\max\{k\geq 0\ |\ [-c_1(X)]^k\neq 0\ \text{in }H^{2k}(X,\mathbb{R})\},\quad 0\leq m\leq n,
\end{equation}
since
\begin{equation}\label{voll}
\frac{1}{n!}\int_X([\omega_0]-2\pi tc_1(X))^n=t^m \frac{(2\pi)^m\binom{n}{m}}{n!}\underbrace{\int_X \omega_0^{n-m}\wedge (-c_1(X))^m}_{>0}+O(t^{m-1}).
\end{equation}
We will now consider separately the behavior of the flow according to the value of the numerical dimension $m$.

\section{The case $m=0$}
The first case to consider is when $m=0$, which by definition means that $c_1(X)=0$ in $H^2(X,\mathbb{R})$. Compact K\"ahler manifolds with this property are known as {\em Calabi-Yau}, and a fundamental theorem of Yau \cite{Ya} shows that every K\"ahler class contains precisely one Ricci-flat K\"ahler metric. In this case, the behavior of the K\"ahler-Ricci flow is well-understood thanks to a classic result of Cao \cite{Ca}, see also \cite[\S 3.4.2--3.4.3]{SWL}, \cite[\S 4.2]{We} for expositions:
\begin{theorem}
Let $X^n$ be a compact Calabi-Yau manifold, and $\omega_0$ a K\"ahler metric on $X$. Then the solution $\omega(t)$ of the K\"ahler-Ricci flow \eqref{krf} starting at $\omega_0$ converges smoothly as $t\to+\infty$ to the unique Ricci-flat K\"ahler metric cohomologous to $\omega_0$.
\end{theorem}
Furthermore, this convergence is exponentially fast (in all $C^k$ norms). A proof of this folklore statement can be found in \cite[p.2941]{TZ}, using the method of \cite{PS}.

\section{The case $m=n$}
Next, we assume that $m=n$. From \eqref{voll} we see that in this case $\vol(X,\omega(t))$ grows to infinity as $t\to+\infty$, and so if we hope to obtain convergence to a reasonable limit we will have to renormalize the flow.  The standard renormalization is to consider
\begin{equation}
\ti{\omega}(s):=e^{-s}\omega(e^s-1),\quad s\in [0,+\infty),
\end{equation}
which now satisfies the K\"ahler-Ricci flow
\begin{equation}\label{krf3}
\frac{\de}{\de s}\ti{\omega}(s)=-\Ric(\ti{\omega}(s))-\ti{\omega}(s),\quad\ti{\omega}(0)=\omega_0.
\end{equation}
and
\begin{equation}
\vol(X,\ti{\omega}(s))=\frac{(2\pi)^n}{n!}\int_X (-c_1(X))^n+O(e^{-s}).
\end{equation}
The condition that $m=n$ is equivalent to $\int_X (-c_1(X))^n>0$, which is usually stated by saying that $K_X$ is nef and big. An important class of manifolds that satisfy this condition are those for which $-c_1(X)\in\mathcal{C}_X$, in which case we say that $K_X$ is ample.

The following result was first obtained by Cao \cite{Ca} when $K_X$ is ample (see  \cite[\S 3.4.1]{SWL}, \cite[\S 4.1]{We} for expositions), and by Tsuji \cite{Ts} and Tian-Zhang \cite{TiZ} in general (see \cite[Theorem 5.3]{To} for an exposition):

\begin{theorem}
Let $X^n$ be a compact K\"ahler manifold with $K_X$ nef and big, and let
\begin{equation}
Z:=\bigcup_{V\subset X,\ \int_V (-c_1(X))^{\dim V}=0}V.
\end{equation}
Then $Z$ is a proper closed analytic subvariety of $X$, and there is a K\"ahler-Einstein metric $\omega_\infty$ on $X\backslash Z$ that satisfies
$$\Ric(\omega_\infty)=-\omega_\infty.$$
Given any  K\"ahler metric $\omega_0$ on $X$, the solution $\ti{\omega}(s)$ of the K\"ahler-Ricci flow \eqref{krf3} starting at $\omega_0$ converges to $\omega_\infty$ locally smoothly on $X\backslash Z$ as $s\to+\infty$.
\end{theorem}
It was also shown by Zhang \cite{Zh2} that the scalar curvature of $\ti{\omega}(s)$ is uniformly bounded independent of $s$.
The subvariety $Z$ is empty if and only if $K_X$ is ample. In this case, the K\"ahler-Einstein metric $\omega_\infty$ was constructed by Aubin \cite{Au} and Yau \cite{Ya}.

In general, a compact K\"ahler manifold with $K_X$ nef and big must be projective (by a classical result of Moishezon, see e.g. \cite[Theorem 2.2.26]{MM}), and by a result of Kawamata-Shokurov \cite[Theorem 2.6]{Ka} there is a birational morphism $f:X\to Y$ onto a normal projective variety, with exceptional locus $\mathrm{Exc}(f)=Z$. It is then shown by Wang \cite{Wa} (see also \cite{GSW, TiZL} for earlier results, and the very recent work of Jian-Song \cite{JS2} for a new proof) that $(X,\ti{\omega}(s))$ converges as $s\to+\infty$ in the Gromov-Hausdorff topology to the metric completion of $(X\backslash Z,\omega_\infty)$, which is a compact metric space homeomorphic to $Y$ by \cite{So}.

\section{The case $0<m<n$: Abundance conjecture}
In order to answer Question \ref{mq}, it thus remains to address the case when $0<m<n$, which is usually referred to as ``intermediate Kodaira dimension''. This is by far the hardest case, and it is open in general. However, decisive progress has been made under the additional assumption that $K_X$ be semiample, which conjecturally always holds, as we now explain.

More precisely, we say that $K_X$ is semiample if there is some $\ell\geq 1$ such that $K_X^\ell$ is globally generated, which means that given any $x\in X$ there is a global section $s\in H^0(X,K_X^\ell)$ with $s(x)\neq 0$. We refer the reader to \cite[\S 2.1.B]{Laz} for a thorough discussion of semiample line bundles, including the basic fact (which follows from \cite[Theorem 2.1.27]{Laz}) that semiample line bundles are always nef.

The fundamental {\em Abundance conjecture} in algebraic geometry (and its natural extension to K\"ahler manifolds) predicts that for the canonical bundle the converse holds as well:
\begin{conjecture}\label{ab}Let $X$ be a compact K\"ahler manifold with $K_X$ nef. Then $K_X$ is semiample.
\end{conjecture}
The Abundance conjecture is known when $n\leq 3$ by \cite{CHP,DO,DO2}.

Motivated by the Abundance conjecture, we will assume in the following section that $K_X$ is not just nef but actually semiample, and discuss recent progress which gives a rather satisfactory picture of the behavior of $(X,\omega(t))$. We will then return to the general case of $K_X$ nef in section \ref{sec}, to discuss what is known and what is expected in that case.

\section{The case $0<m<n$: semiample canonical bundle}\label{semia}
In this section we assume that $X^n$ is a compact K\"ahler manifold with $K_X$ semiample (hence in particular nef). As explained in \cite[Theorem 2.1.27]{Laz}, global sections of $K_X^\ell$ for $\ell$ sufficiently divisible give a holomorphic map
$f:X\to Y,$ with connected fibers onto an irreducible normal projective variety $Y$ with $\dim Y=m$, and $-c_1(X)=f^*[\omega_Y]$ for some K\"ahler metric $\omega_Y$ on $Y$ (understood in the sense of analytic spaces \cite{Mo} if $Y$ is not smooth). In particular,
\begin{equation}
\chi:=f^*\omega_Y,
\end{equation}
is a smooth closed real $(1,1)$-form on $X$ which is semipositive definite and cohomologous to $-2\pi c_1(X)$ (smoothness follows from standard properties of differential forms on analytic spaces).

There is a proper closed analytic subvariety $D\subset Y$ such that $Y^\circ:=Y\backslash D$ is smooth and does not contain any critical values of $f$. This means that if we define $X^\circ:=X\backslash f^{-1}(D)$, then $f:X^\circ\to Y^\circ$ is a proper holomorphic submersion (hence in particular a $C^\infty$ fiber bundle, by Ehresmann's Lemma). The fibers $X_y=f^{-1}(y),y\in Y^\circ$ are automatically Calabi-Yau $(n-m)$-folds (by adjunction) which are pairwise diffeomorphic, but may not be pairwise biholomorphic in general. The fibers of $f$ contained in $f^{-1}(D)$ on the other hand are referred to as singular fibers.

The simplest example of such fiber spaces are products $X=Y\times F$ where $Y$ is a compact K\"ahler manifold with $K_Y$ ample, and $F$ is Calabi-Yau, with the map $f$ being the projection onto the $Y$ factor.

As before, it is most convenient to study the K\"ahler-Ricci flow on $X$ with the normalization as above (now switching notation though)
\begin{equation}\label{krf4}
\frac{\de}{\de t}\omega(t)=-\Ric(\omega(t))-\omega(t),\quad\omega(0)=\omega_0.
\end{equation}
The solution $\omega(t)$ exists for all $t\geq 0$, and is cohomologous to
\begin{equation}
\hat{\omega}(t):=e^{-t}\omega_0+(1-e^{-t})\chi,
\end{equation}
which are K\"ahler metrics on $X$ for all $t\geq 0$.

The study of the K\"ahler-Ricci flow in this setting was initiated by Song-Tian in \cite{ST} when $n=2,m=1$, and in \cite{ST2} in general. Before we can describe the picture in more detail, we need some definitions.

\subsection{The Weil-Petersson form} First, on $Y^\circ$ there is a closed real $(1,1)$-form $\omega_{\rm WP}$ (the {\em Weil-Petersson form}), which is semipositive definite and encodes the variation of the complex structure of the Calabi-Yau fibers $X_y$. In particular, $\omega_{\rm WP}$ is identically zero if and only if all fibers $X_y,y\in Y^\circ$ are biholomorphic to each other. By the Fischer-Grauert theorem \cite{FG}, this happens if and only if $f:X^\circ\to Y^\circ$ is a holomorphic fiber bundle, and in this case we say that $f:X\to Y$ is {\em isotrivial}. We refer the reader to \cite{ST2} or \cite[\S 5.6]{To} for the precise construction of $\omega_{\rm WP}$, which originates in \cite{FS}.

\subsection{The twisted K\"ahler-Einstein metric} Next, by solving an elliptic complex Monge-Amp\`ere equation on $Y$, Song-Tian \cite{ST2} show that there exists a K\"ahler metric $\omega_{\rm can}$ on $Y^\circ$, of the form $\omega_{\rm can}=\omega_Y+\ddbar u$, that satisfies the twisted K\"ahler-Einstein equation
\begin{equation}
  \Ric(\omega_{\rm can})=-\omega_{\rm can}+\omega_{\rm WP}.
\end{equation}
For example, in the above-mentioned product setting where $X=Y\times F$, we have of course $Y^\circ=Y$ and $\omega_{\rm WP}=0$ (since all the fibers are biholomorphic to $F$), so in this case $\omega_{\rm can}$ is simply the K\"ahler-Einstein metric on $Y$ constructed by Aubin and Yau.

\subsection{The semi-Ricci-flat form}The last object that we will need is the semi-Ricci-flat form $\omega_{\rm SRF}$ on $X^\circ$. To define this, for each given $y\in Y^\circ$ we apply the Calabi-Yau theorem to $(X_y, \omega_0|_{X_y})$, which gives us a unique Ricci-flat K\"ahler metric on $X_y$ of the form
\begin{equation}
\omega_0|_{X_y}+\ddbar\rho_y>0, \quad  \int_{X_y}\rho_y\omega_0^{n-m}=0,
\end{equation}
and the smooth function $\rho_y:X_y\to\mathbb{R}$ is uniquely determined thanks to its integral normalization. From Yau's a priori estimates \cite{Ya}, it follows that $\rho_y$ varies smoothly in $y\in Y^\circ$, so letting $y$ vary these define a smooth function $\rho:X^\circ\to \mathbb{R}$, so that $\rho|_{X_y}=\rho_y$. We then define a closed real $(1,1)$-form $\omega_{\rm SRF}$ on $X^\circ$ by
\begin{equation}
  \omega_{\rm SRF}:=\omega_0+\ddbar\rho.
\end{equation}
It is called semi-Ricci-flat because its fiberwise restrictions $\omega_{\rm SRF}|_{X_y}$ are Ricci-flat K\"ahler metrics. However, it is important to note that $\omega_{\rm SRF}$ may fail to be semipositive definite on $X^\circ$, see \cite{CGPT} for a counterexample.

\subsection{Collapsing of the K\"ahler-Ricci flow}With these notations set up, we can now discuss what is known about the behavior of the evolving metrics $\omega(t)$ solving the normalized flow \eqref{krf4} on $X$. The first result, proved by Song-Tian \cite{ST2} is the following:

\begin{theorem}\label{kk}
Let $(X^n,\omega_0)$ be a compact K\"ahler manifold with $K_X$ semiample and numerical dimension $0<m<n$, and let $\omega(t)$ be the solution of \eqref{krf4}.
Then, as $t\to+\infty$, the metrics $\omega(t)$ converge to $f^*\omega_{\rm can}$ in the weak topology of currents on $X$.
\end{theorem}

They also posed the following in \cite[p.612]{ST}, \cite[p.306]{ST2}, \cite[Conjecture 4.5.7]{Ti}, \cite[p.258 and Conjecture 4.7]{Ti2}:

\begin{conjecture}\label{kkk} In the same setting as Theorem \ref{kk}, we have
\begin{itemize}\item[(a)] $\omega(t)\to f^*\omega_{\rm can}$ in the locally smooth topology on $X^\circ$
\item[(b)] Given $K\Subset X^\circ$ there is $C>0$ such that for all $t\geq 0$ we have
\begin{equation}\label{gorm}
\sup_K|\Ric(\omega(t))|_{\omega(t)}\leq C.
\end{equation}
\end{itemize}
\end{conjecture}
To clarify, item (a) means that given any $K\Subset X^\circ$  and $k\geq 0$, we have
\begin{equation}\label{germ}
\|\omega(t)-f^*\omega_{\rm can}\|_{C^k(K,\omega_0)}\to 0, \quad \text{as }t\to+\infty.
\end{equation}
while item (b) simply means that the Ricci curvature of $\omega(t)$ remains uniformly bounded on $K$, independent of $t$.

The first progress towards establishing (a) was done by Fong-Zhang \cite{FZ}, adapting a method by the author \cite{To2} in a related elliptic PDE, who showed that the K\"ahler potentials of $\omega(t)$ (with respect to $\hat{\omega}(t)$) converge to $f^*u$ (the potential of $f^*\omega_{\rm can})$ in $C^{1,\alpha}(K,\omega_0)$ for all $\alpha<1$. This however falls short of proving \eqref{germ} for $k=0$, and this was only achieved later by Weinkove, Yang and the author in \cite{TWY} with substantial more work (see \cite[\S 5.5--5.13]{To} for a detailed exposition). Using this result, Zhang and the author \cite{TZ} identified the next order behavior of $\omega(t)$ when restricted to a smooth fiber, which is
\begin{equation}
\|(\omega(t)-e^{-t}\omega_{\rm SRF})|_{X_y}\|_{C^k(X_y,\omega_0|_{X_y})}=o(e^{-t}),
\end{equation}
for any $y\in Y^\circ$ and $k\geq 0$. Next, it was observed in \cite{FZ,HT,TZ} (see also \cite[\S 5.14]{To} for a unified exposition, and \cite{Gi} for the case of products) that Conjecture \ref{kkk} (a) holds in the case when the smooth fibers $X_y$ are tori or finite quotients of tori. This used crucially ideas of Gross, Zhang and the author \cite{GTZ} in the elliptic setting. Further progress was made by Fong-Lee \cite{FL} who proved Conjecture \ref{kkk} (a) when $f$ is isotrivial, and by Chu-Lee \cite{CL} who proved \eqref{germ} for the H\"older norm $C^\alpha$, $0<\alpha<1$. Both of these works are based on ideas of Hein and the author \cite{HT2} in the elliptic case.

About Conjecture \ref{kkk} (b), Song-Tian \cite{ST3} showed that the scalar curvature of $\omega(t)$ is uniformly bounded on all of $X$, while Zhang and the author \cite{TZ} showed that the full Riemann curvature tensor of $\omega(t)$ remains uniformly bounded on compact subsets $K\Subset X^\circ$ if and only if the smooth fibers $X_y$ are tori or finite quotients of tori (see also \cite[p.364]{To}). Thus, the conjectured Ricci bound \eqref{gorm} is in some sense optimal. The bound \eqref{gorm} was established in \cite{FL} when $f$ is isotrivial.

Very recently, Hein, Lee and the author \cite{HLT} were able to confirm Conjecture \ref{kkk} in general:

\begin{theorem}Conjecture \ref{kkk} holds.
\end{theorem}

This result can be thought of as a parabolic version of the work \cite{HT3} by Hein and the author, in the elliptic setting. To prove \eqref{germ}, the authors prove a much more precise asymptotic expansion for $\omega(t)$, as follows. Given $k\geq 0$, and given a sufficiently small coordinate ball $B\subset Y^\circ$ (over which $f$ is smoothly trivial, so $f^{-1}(B)$ is diffeomorphic to $B\times F$ for some Calabi-Yau manifold $F$), they show that, up to shrinking $B$, we can write on $B\times F$ for all $t\geq 0$,
\begin{equation}\label{exp}
\omega(t)=f^*\omega_{\rm can}+e^{-t}\omega_{\rm SRF}+\gamma_0(t)+\sum_{j=1}^k\gamma_{j,k}(t)+\eta_k(t),
\end{equation}
where $\gamma_0(t)$ is pulled back from $B$ and goes to zero in $C^k(B)$, the $\gamma_{j,k}(t),1\leq j\leq k$ go to zero in $C^k(B\times F, \omega_0)$ (at different speeds that depend on $j$), and the remainder $\eta_k(t)$ goes to zero in an even stronger $t$-dependent ``shrinking'' $C^k$ norm on $B\times F$. This description clearly implies Conjecture \ref{kkk} (a), and then plugging in this expansion with $k=2$ into the definition of $\Ric(\omega(t))$ and employing explicit estimates for the pieces of the decomposition \eqref{exp} we show Conjecture \ref{kkk} (b).

In the two aforementioned special cases when $f$ is isotrivial, or when the smooth fibers are tori or finite quotients, it turns out that all the terms $\gamma_{j,k}(t)$ in \eqref{exp} vanish identically, so the result is stronger (and the proof much easier). In general however these terms do not vanish, and simply identifying them requires very substantial work. Once the pieces of the decomposition \eqref{exp} have been defined, the required {\em a priori} estimates on all the pieces are proved by an intricate argument by contradiction and blowup, in the spirit of \cite{HT3}. Furthermore, the estimates proved are actually in stronger parabolic H\"older norms.

\subsection{Gromov-Hausdorff limits}
While the results in the previous subsection give very strong estimates for $\omega(t)$ away from the singular fibers of $f$, these estimates blow up much too fast near the singular fibers to give us useful information. Nevertheless, Song-Tian \cite{ST,ST2,ST4} conjectured that $(X,\omega(t))$ should have a Gromov-Hausdorff limit as $t\to+\infty$, as follows:

\begin{conjecture}\label{lll}
In the same setting as Theorem \ref{kk}, we have
\begin{itemize}\item[(a)] There is $C>0$ such that for all $t\geq 0$,
\begin{equation}\label{dia}\mathrm{diam}(X,\omega(t))\leq C,\end{equation}
\item[(b)] As $t\to+\infty$,  $(X,\omega(t))$ converges in the Gromov-Hausdorff topology to $(Z,d)$, the metric completion of $(Y^\circ,\omega_{\rm can})$
\item[(c)] $Z$ is homeomorphic to $Y$.
\end{itemize}
\end{conjecture}

After much work on these questions, Conjecture \ref{lll} (a) was completely settled by Jian-Song \cite{JS}. As a consequence, they also deduced that sequential Gromov-Hausdorff limits exist, which is not a priori clear. Another breakthrough was achieved by Song-Tian-Zhang \cite{STZ}, who proved Conjecture \ref{lll} (b) and (c) in the case when $m=1$ and $X_y$ tori, and combining their method with \cite{JS} settles Conjecture \ref{lll} (b) and (c) completely when $m=1$. Furthermore, \cite{STZ} proved Conjecture \ref{lll} (c) when $Y$ is smooth or has orbifold singularities. Lastly,  Li and the author \cite{LT} proved Conjecture \ref{lll} (b) for arbitrary $m$ assuming that $Y$ is smooth and that the union $D^{(1)}$ of the codimension $1$ irreducible components of $D$ is a simple normal crossings divisor (recall that $D$ was defined at the beginning of Section \ref{semia}).  This last assumption (which always holds when $m=1$) is used to appeal to the estimates in \cite{GTZ2} which show that in this case $\omega_{\rm can}$ has conical singularities along $D^{(1)}$ (up to a small logarithmic error), and it extends smoothly across $D\backslash D^{(1)}$.

\section{The case $0<m<n$: nef canonical bundle}\label{sec}
In this section we assume that $(X^n,\omega_0)$ is a compact K\"ahler manifold with $K_X$ nef and numerical dimension $0<m<n$, but we do not assume that it is semiample (although it is predicted to be so by the Abundance Conjecture \ref{ab}). In this situation much less is known about the long-time behavior of the solution $\omega(t)$ of the K\"ahler-Ricci flow \eqref{krf4}.

\subsection{Weak limit of the flow}
To start, we now let $\chi:=-\Ric(\omega_0)$, which is a closed real $(1,1)$-form cohomologous to $-2\pi c_1(X)$, and let
\begin{equation}\alpha(t)=e^{-t}\omega_0+(1-e^{-t})\chi,
\end{equation}
which is a family of closed real $(1,1)$-forms cohomologous to $\omega(t)$, with no positivity property in general. It is again easy to see \cite[(5.30)]{To} that \eqref{krf4} holds if and only if we can write
\begin{equation}
\omega(t)=\alpha(t)+\ddbar\vp(t),
\end{equation}
with $\vp(t)$ solving the parabolic Monge-Amp\`ere equation
\begin{equation}\label{ma}
\left\{
                \begin{aligned}
                  &\frac{\de}{\de t}\vp(t)=\log\frac{e^{(n-m)t}(\alpha(t)+\ddbar\vp(t))^n}{\omega_0^n}-\vp(t)\\
                  &\vp(0)=0\\
                  &\alpha(t)+\ddbar\vp(t)>0.
                \end{aligned}
              \right.
\end{equation}
on $X\times [0,+\infty)$. Let also $h$ be the Hermitian metric on $K_X$ naturally induced by $\omega_0$, whose curvature form is equal to $\chi$.

The following result is a simple consequence of work of Fu-Guo-Song \cite{FGS} (see also \cite{GPSS,GPTW,GS}):
\begin{theorem}\label{known}
Let $(X^n,\omega_0)$ be a compact K\"ahler manifold with $K_X$ nef  and numerical dimension $0<m<n$, and let $\omega(t)$ be the immortal solution of \eqref{krf4}. Then
\item[(a)] Given any sequence $t_i\to\infty$ we can find a subsequence and a quasi-psh function $\vp_\infty:X\to\mathbb{R}\cup\{-\infty\}$, with $\chi_\infty:=\chi+\ddbar\vp_\infty\geq 0$ weakly, such that $\vp(t_i)\to\vp_\infty$ in $L^1(X)$ as $t\to\infty$. In particular $\omega(t_i)\to \chi_\infty$ weakly as currents. A priori $\chi_\infty$ may depend on the sequence.
\item[(b)] $\chi_\infty$ is a closed positive $(1,1)$-current with minimal singularities in the class $-2\pi c_1(X)$
\item[(c)] Given any $\ell\geq 1$ and any $s\in H^0(X,K_X^\ell)$, we have
\begin{equation}\label{l2}
\int_X |s|^2_{h^\ell}e^{-\ell\vp_\infty}\omega_0^n<\infty.
\end{equation}
\end{theorem}

In part (b) the notion of minimal singularities in the class $[\chi]=-2\pi c_1(X)$, introduced in \cite{DPS}, means that given any quasi-psh function $\eta:X\to\mathbb{R}\cup\{-\infty\}$ with $\chi+\ddbar\eta\geq 0$ weakly, there is $C>0$ such that
\begin{equation}
\vp_\infty\geq \eta -C,
\end{equation}
holds on $X$.

Property (c) is called ``analytic Zariski decomposition''  by Tsuji \cite{Ts}. Note that a priori it may happen that $H^0(X,K_X^\ell)=0$ for all $\ell\geq 1$, although this is predicted not to happen as a consequence of the Abundance Conjecture.

\begin{proof}
(a) It is well-known (see e.g. \cite[Lemma 4.3]{GS}) that there is $C>0$ such that for all $t\geq 0$,
\begin{equation}\label{scr}
  \sup_X \vp(t)\leq C, \quad \sup_X\dot{\vp}(t)\leq C.
\end{equation}
Thanks to the uniform upper bound for $\vp(t)$, the statement in part (a) follows from basic compactness properties of quasi-psh functions (see e.g. \cite[Theorem 3.2.12]{Hor}), provided we show that $\vp$ does not converge to $-\infty$ uniformly on $X$. To see this, we use the Monge-Amp\`ere equation \eqref{ma} to get
\begin{equation}
\int_X e^{\vp+\dot{\vp}}\omega_0^n=e^{(n-m)t}\int_X\omega(t)^n=n!e^{(n-m)t}\vol(X,\omega(t)),
\end{equation}
and using Stokes we have
\begin{equation}\begin{split}
\vol(X,\omega(t))&=\frac{1}{n!}\int_X(e^{-t}\omega_0-2\pi (1-e^{-t})c_1(X))^n\\
&=e^{-(n-m)t}\frac{(2\pi)^m\binom{n}{m}}{n!}\underbrace{\int_X \omega_0^{n-m}\wedge (-c_1(X))^m}_{>0}+O(e^{-(n-m+1)t}),
\end{split}
\end{equation}
and so
\begin{equation}
\int_X e^{\vp+\dot{\vp}}\omega_0^n\geq C^{-1},
\end{equation}
which implies
\begin{equation}\label{klk}
\sup_X (\vp(t)+\dot{\vp}(t))\geq -C.
\end{equation}
But since $\dot{\vp}(t)\leq C$, then \eqref{klk} implies that $\sup_X\vp(t)\geq -C$, as desired.

(b) By part (a) we know that given $t_i\to\infty$ we can find a subsequence and $\vp_\infty$ quasi-psh with $\chi+\ddbar\vp_\infty\geq 0$ such that $\vp(t_i)\to \vp_\infty$ in $L^1(X)$. We need to show that $\vp_\infty$ has minimal singularities. For this, we consider the envelope
$$V_\infty(x)=\sup\{\psi(x)\ |\ \chi+\ddbar\psi\geq 0, \psi\leq 0\},$$
which satisfies $V_\infty\leq 0$, is quasi-psh with $\chi+\ddbar V_\infty\geq 0$, and has minimal singularities: given any quasi-psh function $\eta:X\to\mathbb{R}\cup\{-\infty\}$ with $\chi+\ddbar\eta\geq 0$ weakly, since $\eta$ is usc there is $C>0$ such that $\eta\leq C$ on $X$, so taking $\psi=\eta-C$ in the definition of the envelope shows that $V_\infty\geq \eta-C,$ as desired. We empasize here that we do not know whether $V_\infty$ is bounded on $X$.
If we show that
\begin{equation}\label{li}
\|\vp_\infty-V_\infty\|_{L^\infty(X)}\leq C,
\end{equation}
then it follows immediately that $\vp_\infty$ also has minimal singularities.
To prove \eqref{li}, we first use \eqref{scr} and see that
\begin{equation}
\omega(t)^n\leq Ce^{-(n-m)t}\omega_0^n.
\end{equation}
We can then appeal to \cite[Prop. 1.1]{FGS} (which was recently reproved and improved in \cite[Theorem 1]{GPTW}) and obtain that
\begin{equation}\label{linf}
\|\vp(t)-V_t\|_{L^\infty(X)}\leq C,
\end{equation}
for all $t\geq 0$,
where
\begin{equation}
V_t(x)=\sup\{\psi(x)\ |\ \alpha(t)+\ddbar\psi\geq 0, \psi\leq 0\},
\end{equation}
which satisfies $V_t\leq 0$, is quasi-psh, and $\alpha(t)+\ddbar V_t\geq 0$ weakly. Since the class $\alpha(t)$ is K\"ahler, the functions $V_t$ are in $C^{1,1}(X)$ by \cite{CZ,To3}.

We now claim that as $t\to\infty,$ the functions $V_t$ converge pointwise to $V_\infty$.
Indeed, we clearly have that $(1-e^{-t})V_\infty$ participates to the supremum defining $V_t$, and hence
\begin{equation}\label{kuz}
V_\infty\leq \frac{V_t}{1-e^{-t}},
\end{equation}
while on the other hand for $t<s$, if $\psi\leq 0$ is quasi-psh and satisfies $\alpha(s)+\ddbar\psi\geq 0$ weakly, then we also have
\begin{equation}
\alpha(t)+\frac{1-e^{-t}}{1-e^{-s}}\ddbar\psi\geq \frac{1-e^{-t}}{1-e^{-s}}(\alpha(s)+\ddbar\psi)\geq 0,
\end{equation}
using that $\frac{e^{-s}}{1-e^{-s}}\leq \frac{e^{-t}}{1-e^{-t}}$,
and so $\frac{1-e^{-t}}{1-e^{-s}}\psi$ participates to the supremum defining $V_t$, hence
\begin{equation}
\frac{\psi}{1-e^{-s}}\leq  \frac{V_t}{1-e^{-t}},
\end{equation}
and taking here the supremum over all such $\psi$ gives
\begin{equation}
\frac{V_s}{1-e^{-s}}\leq  \frac{V_t}{1-e^{-t}}.
\end{equation}
This shows that $\frac{V_t}{1-e^{-t}}$ decrease pointwise as $t\to\infty$ to some limit usc function $V'_\infty:X\to [-\infty,+\infty)$ which by \eqref{kuz} satisfies
\begin{equation}\label{cuz}
V_\infty\leq V'_\infty\leq 0,
\end{equation}
so in particular it is not identically $-\infty$. Observe that by dominated convergence we have $\frac{V_t}{1-e^{-t}}\to V'_\infty$ in $L^1(X)$, and so $\frac{1}{1-e^{-t}}\ddbar V_t\to \ddbar V'_\infty$ weakly as currents.  Since
\begin{equation}
\frac{1}{1-e^{-t}}(\alpha(t)+\ddbar V_t)\geq 0,
\end{equation}
weakly, passing this to the weak limit as $t\to+\infty$ shows that $\chi+\ddbar V'_\infty\geq 0$ weakly, so $V'_\infty$ participates to the supremum defining $V_\infty$, and hence
\begin{equation}
V'_\infty\leq V_\infty,
\end{equation}
and from this and \eqref{cuz} we conclude that $V'_\infty=V_\infty$. We have thus shown that $\frac{V_t}{1-e^{-t}}\to V_\infty$ pointwise as $t\to\infty$, and hence also $V_t\to V_\infty$ pointwise, as claimed.

Passing to our sequence $t_i\to\infty$ as in part (a), we have that $\vp(t_i)$ converges to $\vp_\infty$ in $L^1(X)$, and (up to passing to a further subsequence) also pointwise a.e. Passing to the limit in \eqref{linf}, and using the above claim, gives
\begin{equation}\label{linfi}
|\vp_\infty-V_\infty|\leq C,
\end{equation}
a.e. on $X$, hence everywhere by elementary properties of psh functions (see e.g. \cite[Theorem K.15]{Gu}), and \eqref{li} follows.

(c) Given any $\ell\geq 1$ and any nontrivial $s\in H^0(X,K_X^\ell)$ (if it exists), since the metric $h^\ell$ on $K_X^\ell$ has curvature equal to $\ell\chi$, it follows from the Poincar\'e-Lelong equation that
\begin{equation}
\chi+\frac{1}{\ell}\ddbar\log |s|^2_{h^\ell}\geq 0,
\end{equation}
weakly. Thus, since $\vp_\infty$ has minimal singularities, it follows that there exists $C>0$ such that
\begin{equation}
\vp_\infty\geq \frac{1}{\ell}\log |s|^2_{h^\ell}-C,
\end{equation}
or equivalently,
\begin{equation}
|s|^2_{h^\ell}e^{-\ell\vp_\infty}\leq C,
\end{equation}
and integrating this against $\omega_0^n$ gives \eqref{l2}.
\end{proof}

The following is then a very natural expectation (see also \cite{Zh}):

\begin{conjecture}\label{unknown}
Let $(X^n,\omega_0)$ be a compact K\"ahler manifold with $K_X$ nef  and numerical dimension $0<m<n$, and let $\omega(t)$ be the immortal solution of \eqref{krf4}. Then there is a quasi-psh function $\vp_\infty$ with $\chi_\infty:=\chi+\ddbar\vp_\infty\geq 0$ weakly, such that $\vp(t)\to\vp_\infty$ in $L^1(X)$ as $t\to+\infty$. In particular, $\omega(t)\to \chi_\infty$ weakly as currents as $t\to+\infty$. Furthermore, $\chi_\infty$ is independent also of the choice of the initial metric $\omega_0$.
\end{conjecture}
In other words, the subsequential limits $\vp_\infty$ in Theorem \ref{known} should be unique, and independent of the subsequence, and the current $\chi_\infty$ should be even independent of the initial metric. Of course, by Theorem \ref{known}, such $\chi_\infty$ would have minimal singularities.  If Conjecture \ref{unknown} is settled, then $he^{-\vp_\infty}$ would be a singular metric on $K_X$ with semipositive curvature current $\chi_\infty$, which up to scaling by a constant would be completely canonical, i.e. depending only on the complex structure of $X$.

Such a singular metric could be used for $L^2$ type arguments, especially if one can establish further regularity of $\vp_\infty$. When $K_X$ is semiample, it is known by \cite{DZ} that $\vp_\infty$ is continuous on $X$, and it may even be H\"older continuous. One can then ask about the following weaker properties:

\begin{conjecture}
In the same setting as Conjecture \ref{unknown}, the following hold:
\begin{itemize}
\item[(a)] $\vp_\infty$ has vanishing Lelong number at every point
\item[(b)] $\vp_\infty$ is bounded
\end{itemize}
\end{conjecture}
By a result of Skoda \cite{Sk}, if (a) holds then $e^{-A\vp_\infty}\in L^1(X)$ for all $A>0$, which implies that \eqref{l2} holds. Thus, condition (a) can be viewed as a strengthening of \eqref{l2}.
The implication (b)$\Rightarrow$(a) is well-known, and as mentioned above (b) holds when $K_X$ is semiample.

These conditions can also be stated in terms of the smooth potentials $\vp(t)$. Thanks to a result of Demailly-Koll\'ar \cite{DK}, condition (a) is equivalent to the fact that for every $A>0$ there is $C>0$ such that for all $t\geq 0$,
\begin{equation}
\int_X e^{-A\vp(t)}\omega_0^n\leq C,
\end{equation}
while condition (b) is equivalent to the existence of $C>0$ such that for all $t\geq 0$,
\begin{equation}
\inf_X\vp(t)\geq -C.
\end{equation}

\subsection{Singularity type at infinity}
Following Hamilton \cite{Ha2}, we say that an immortal solution of \eqref{krf3} develops a ``Type III'' singularity at infinity if there is $C>0$ such that for all $t\geq 0$ we have
\begin{equation}
\sup_X|\mathrm{Rm}(\omega(t))|_{\omega(t)}\leq C,
\end{equation}
and develops a ``Type IIb'' singularity if this fails. In \cite{TZ}, Zhang and the author gave an almost complete classification (complete when $n=2$) of what the possible singularity type is, under the assumption that $K_X$ is semiample, see also \cite[Theorem 6.6]{To}. In all the cases that they considered, the singularity type was actually independent of the initial metric $\omega_0$, and the author conjectured in \cite[Conjecture 6.7]{To} that this should always be the case even just assuming that $K_X$ is nef. Zhang \cite{ZhY2} proved this conjecture when $K_X$ is semiample, and very recently Wondo-Zhang \cite{WZ} settled it completely.

It remains a very interesting problem to complete the classification of \cite{TZ} in all dimensions, under the assumption that $K_X$ is semiample. The only case which is left to understand is when the fiber space $f$ has some singular fibers, and the smooth fibers $X_y$ are all tori or finite quotients of tori.

\subsection{Diameter bounds}
Lastly, let us mention a recent striking result of Guo-Phong-Song-Sturm \cite{GPSS}:
let $(X^n,\omega_0)$ be a compact K\"ahler manifold with $K_X$ nef  and numerical dimension $0<m<n$, and let $\omega(t)$ be the immortal solution of \eqref{krf4}. If we also suppose that $H^0(X,K_X^\ell)\neq 0$ for some $\ell\geq 1$ (i.e. that $X$ has nonnegative Kodaira dimension), then there is $C>0$ such that for all $t\geq 0$ we have
\begin{equation}
\mathrm{diam}(X,\omega(t))\leq C,
\end{equation}
and furthermore, given any $t_i\to+\infty$, up to passing to a subsequence we will have that $(X,\omega(t_i))$ converges in Gromov-Hausdorff to some compact metric space $(Z,d)$.  It would be of course very interesting to remove the condition that $X$ have nonnegative Kodaira dimension (which is expected to hold by the Abundance conjecture).

\end{document}